\documentclass[final,3p,times]{elsarticle}
\usepackage[latin1]{inputenc}
\usepackage[swedish,english]{babel}
\usepackage{url}
\usepackage{amsmath,amsthm}
\usepackage{amssymb}
\usepackage{graphics}
\usepackage{subfigure}
\usepackage{sidecap}

\numberwithin{equation}{section}
\numberwithin{table}{section}
\numberwithin{figure}{section}

\theoremstyle{plain}

\theoremstyle{definition}
\newtheorem{definition}{Definition}[section]
\newtheorem{assumption}[definition]{Assumption}

\newtheorem{criterion}[definition]{Criterion}

\theoremstyle{remark}

\newcommand{\thetacriterion}{$\theta$-criterion}
\newcommand{\realdom}{\mathbf{R}}
\newcommand{\intdom}{\mathbf{Z}}

\newcommand{\Gfar}{\widehat{G}}
\newcommand{\Glocal}{\widetilde{G}}

\newcommand{\Eff}{E_{\mbox{{\footnotesize far-to-far}}}}
\newcommand{\Efl}{E_{\mbox{{\footnotesize far-to-local}}}}

\newcommand{\Ordo}[1]{\mathcal{O}\left(#1\right)}
\newcommand{\ontop}[2]{\genfrac{}{}{0pt}{}{#1}{#2}}
\newcommand{\TOL}{\operatorname{TOL}}

\biboptions{sort&compress}

\begin{document}

\selectlanguage{english}

\begin{frontmatter}
  
\title{On well-separated sets and fast multipole methods}

\author{Stefan Engblom\corref{corr}}

\cortext[corr]{Phone: +46-18-471 27 54, Fax: +46-18-51 19 25}

\ead{stefane@it.uu.se}
\ead[url]{http:\slash\slash user.it.uu.se/$\sim$stefane} % quite ugly!

%\date{June 27, 2011}

\address{Department of Numerical Analysis, \\
  School of Computer Science and Communication\fnref{change}, \\
  Royal Institute of Technology, S-100 44 Stockholm, Sweden}

% the footnote-text shows up but not the actual footnote
\fntext[change]{Current address: Division of Scientific Computing,
  Department of Information Technology, Uppsala University, Box 337,
  SE-751 05 Uppsala, Sweden.}

\begin{abstract}

  The notion of well-separated sets is crucial in fast multipole
  methods as the main idea is to approximate the interaction between
  such sets via cluster expansions. We revisit the one-parameter
  multipole acceptance criterion in a general setting and derive a
  relative error estimate. This analysis benefits asymmetric versions
  of the method, where the division of the multipole boxes is more
  liberal than in conventional codes. Such variants offer a
  particularly elegant implementation with a balanced multipole tree,
  a feature which might be very favorable on modern computer
  architectures.

\end{abstract}


% 65-xx Numerical analysis
% 65Mxx Partial differential equations, initial value and
%       time-dependent initial-boundary value problems
% 65M15 Error bounds
% 65M80 Fundamental solutions, Green's function methods, etc.

\end{frontmatter}

%**************************************************************************

\section{Motivation}

We consider in this paper a general error analysis and some
implementation issues for fast multipole methods (FMMs). Since their
first appearance in \cite{FMM,AFMM}, these tree-based algorithms have
become important computational tools for evaluating pairwise
interactions of the type
\begin{align}
  \Phi(x_{i}) &= \sum_{j = 1, j \not = i}^{N} G(x_{i},x_{j}),
  \quad x_{i} \in \realdom^{D}, \quad i = 1 \ldots N.
\end{align}
FMMs offer an $\Ordo{N}$-complexity and an \textit{a priori} error
estimate, but implementing a fully fledged adaptive FMM in 3D is a
daunting task \cite{new_AFMM}. Parallelization issues complicate
matter even further and call for a balance between, on the one hand
theoretical efficiency, and on the other hand software complexity
\cite{ace,pace}. A major inconvenience with \emph{adaptive} versions
is the complicated memory access pattern which is due to the
communication between levels in the multipole tree. Although this can
be mitigated either through post-balancing algorithms
\cite{balancedtrees}, or by employing more advanced data-structures
\cite[Chap.~6.6]{helmholtzFMM}, with modern data-parallel
architectures it has in fact been suggested that uniform versions
offer better performance \cite{fmmgpu}.

\begin{SCfigure}[][htb]
  \includegraphics{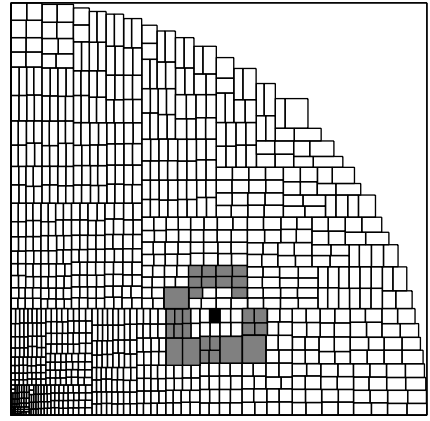}
  \caption{Sample asymmetric adaptive multipole mesh in 2D. The
    discretization is obtained by successively splitting boxes along
    each coordinate axis in such a way that the number of sources in
    the four resulting boxes is very nearly equal. In this example the
    grey boxes are to interact through cluster-to-cluster interactions
    with the black box. That is, they satisfy the
    \thetacriterion\ ($\theta = 1/2$).}
  \label{fig:mesh}
\end{SCfigure}

An alternative is to use \emph{asymmetric} adaptive meshes as outlined
in Figures~\ref{fig:mesh} and \ref{fig:boxes}. Here the tree becomes
balanced at the cost of a variable, but local, communication stencil.
Also, the actual form of the \emph{multipole acceptance criterion}
becomes more critical as the mesh looses regularity. Comparisons
between different criteria for cell-to-particle methods are found in
\cite{fmmskeletons}, while a careful worst-case analysis of uniform
FMMs is found in \cite{fmmerr}.

The main contributions of the current paper are found in
Sections~\ref{sec:assumptions} and \ref{sec:analysis} where the
precise statement of the multipole acceptance criterion and the
required assumptions on the kernel $G$ are discussed together with a
general error analysis. Our approach is inspired by the treatment in
the monograph \cite{griebelMoldyn}, but we are able to offer several
important corrections. In Sections~\ref{sec:implementation} and
\ref{sec:experiments} implementation issues are highlighted and we
also perform numerical experiments illustrating the sharpness of the
theory and the efficiency of the proposed approach.

%**************************************************************************

\section{Well-separated sets and kernel assumptions}
\label{sec:assumptions}

We begin by stating the one-parameter criterion for two sets being
well-separated. This is also a correction to the too weak criterion
found in \cite[Eq.~(8.31), p.~379]{griebelMoldyn}.

\begin{criterion}[\thetacriterion]
  \label{crit:theta}
  Let the sets $S_{1} , S_{2} \subset \realdom^{D}$ be contained
  inside two disjoint spheres such that $\|S_{1}-x_{0}\| \le r_{1}$
  and $\|S_{2}-y_{0}\| \le r_{2}$. Given $\theta \in (0,1)$, if $d
  \equiv \|x_{0}-y_{0}\|$, $R \equiv \max\{r_{1},r_{2}\}$, and $r
  \equiv \min\{ r_{1},r_{2}\}$, then the two sets are
  \emph{well-separated} whenever $R+\theta r \le \theta d$.
\end{criterion}
In other words, any of the two sets may be expanded by a factor of
$1/\theta$ and arbitrarily rotated about its center point without
touching the other set. When regarded as a parameter inside the FMM it
will be evident that a smaller $\theta$ yields a smaller error at the
cost of a larger communication stencil.

We remark that the \thetacriterion\ under consideration is symmetric
in its two arguments as a reflection upon the fact that we mainly
consider FMMs where the sources and the evaluation points are
``roughly'' the same. --- Indeed, our implementation as outlined in
Section~\ref{sec:implementation} produces a representation of the
total field in the \emph{whole} enclosing box under consideration. An
algorithm which is adaptive in sources and evaluation points
separately is described in \cite[Chap.~6.6; see also
Fig.~5.1]{helmholtzFMM}.

We shall need the following two simple consequences of the
\thetacriterion: since $\theta < 1$ we get
\begin{align}
  \label{eq:ineq1}
  \frac{d+r+R}{d-R} &\le \frac{d+r+R/\theta}{d-\theta(d-r)} \le
  \frac{2d}{(1-\theta)d+\theta r} \le \frac{2}{1-\theta}.
  \intertext{Furthermore, writing the \thetacriterion\ as $r \le
    \theta(d-R)+(\theta-1)(R-r)$, we also see that}
  \label{eq:ineq2}
  \frac{r}{d-R} &\le \frac{R}{d-r} \le \theta.
\end{align}

For $\alpha$ and $\beta$ multi-indices in the normed space
$(\intdom_{+}^{D},|\cdot|)$ we write for brevity
$G_{\alpha,\beta}(x,y) = \partial^{\alpha} x \partial^{\beta} y
G(x,y)$ and additionally define factorials and powers in the usual
intuitive manner. Similarly to \cite[Eq.~(8.15),
  p.~317]{griebelMoldyn} (but here with an explicit factor $n!$) we
assume that the kernel $G$'s first few derivatives are equivalent to
the harmonic potential:
\begin{assumption}[Kernel regularity]
  \label{ass:Gbound}
  For $\alpha,\beta \in \intdom_{+}^{D}$ with $|\alpha|, |\beta| \le
  p+1$,
\begin{align}
  \label{eq:dGbound}
  |G_{\alpha,\beta}(x,y)| &\le C \frac{n!}{\|x-y\|^{n+1}}, \\
  \intertext{where $n \equiv |\alpha+\beta|$. Additionally, $G$ is
    positive and satisfies}
  \label{eq:Gbound}
  \|x-y\|^{-1} &\le cG(x,y).
\end{align}
\end{assumption}

The analysis below takes place in $D$-dimensional Cartesian space. In
order to avoid the worst-case bound $\|\cdot\|_{1} \le \sqrt{D}
\|\cdot\|$, assuming some kind of symmetry of the kernel seems
inevitable (see also \cite{ace}):
\begin{assumption}[Rotational invariance]
  \label{ass:Grot}
  For any rotation $T$ of the coordinate system,
  \begin{align}
    \label{eq:Grot}
    |G(x,y)| = |G(Tx,Ty)|.
  \end{align}
\end{assumption}
In other words, we may freely rotate the coordinate system before
expanding the kernel (provided of course that any expansion points are
rotated as well). By choosing a $T$ such that $\|Tx\|_{1} = \|x\|$,
results obtained below in $\|\cdot\|_{1}$ will be transferred to the
Euclidean norm without introducing any constants.

%**************************************************************************

\section{Analysis}
\label{sec:analysis}

Consider now two points $x \in S_{1}$ and $y \in S_{2}$, where $S_{1}$
and $S_{2}$ satisfy the \thetacriterion. Approximating a unit source
at $y$ using a far-to-far translation ($G \to \Gfar_{p}$, centered at
$y_{0}$), followed by a far-to-local expansion ($\Gfar_{p} \to
\Glocal_{p}$, centered at $x_{0}$), can be written as
\begin{align}
  \nonumber
  G(x,y) &= \Glocal_{p}(x,y; \; x_{0},y_{0})-
  \left[ \Glocal_{p}(x,y; \; x_{0},y_{0})-
    \Gfar_{p}(x,y; \; y_{0}) \right]-
  \left[ \Gfar_{p}(x,y; \; y_{0})-G(x,y) \right] \\
  \label{eq:errsplit}
  &=: \Glocal_{p}(x,y; \; x_{0},y_{0})-\Efl-\Eff,
\end{align}
with $p$ the order of the expansion. We consider the two errors in
turn.

The integral form for the remainder of the $D$-dimensional Taylor
series becomes
\begin{align*}
  |\Eff| &= (p+1) \left| \sum_{|\beta| = p+1}
  \frac{(y-y_{0})^{\beta}}{\beta!}
  \int_{0}^{1} (1-t)^{p} G_{0,\beta}(x,y_{0}+t(y-y_{0}))
  \, dt \right| \\
  &\le C (p+1) \frac{\|y-y_{0}\|_{1}^{p+1}}{\|x-y_{0}\|^{p+2}}
  \int_{0}^{1} \frac{(1-t)^{p}}{\left( 1-t
    \frac{\|y-y_{0}\|}{\|x-y_{0}\|} \right)^{p+2}} \, dt,
\end{align*}
using \eqref{eq:dGbound} and the multinomial theorem. As discussed in
conjunction with Assumption~\ref{ass:Grot} above we may now replace
$\|y-y_{0}\|_{1}$ with $\|y-y_{0}\|$. Using the triangle inequality
and \eqref{eq:ineq2} yields $\|y-y_{0}\|/\|x-y_{0}\| \le \theta$ so
that from \eqref{eq:Gbound},
\begin{align}
  \label{eq:relG1}
  \frac{1}{\|x-y_{0}\|} &\le cG(x,y) \frac{\|x-y\|}{\|x-y_{0}\|}
  \le cG(x,y) \left( 1+\frac{\|y-y_{0}\|}{\|x-y_{0}\|} \right)
  \le cG(x,y) (1+\theta).
\end{align}
Using $1-t \le 1-\theta t \le 1-t \, \|y-y_{0}\|/\|x-y_{0}\|$ we
readily bound the integrand and get
\begin{align}
  \label{eq:Ef2f}
  |\Eff| &\le cC (p+1) \; \theta^{p+1} \; 
  \frac{1+\theta}{1-\theta} G(x,y).
\end{align}

As for the first term in \eqref{eq:errsplit} we obtain this time a
\emph{sum} of integral remainders,
\begin{align*}
  |\Efl| &= (p+1) \left| \sum_{\ontop{|\alpha| = p+1}{|\beta| \le p}} 
  \frac{(x-x_{0})^{\alpha}}{\alpha!} \frac{(y-y_{0})^{\beta}}{\beta!}
  \int_{0}^{1} (1-t)^{p} G_{\alpha,\beta}(x_{0}+t(x-x_{0}),y_{0})
  \, dt \right| \\
  &\le C (p+1) \sum_{k = 0}^{p} 
  \frac{\|x-x_{0}\|^{p+1}}{(p+1)!} \frac{\|y-y_{0}\|^{k}}{k!}
  \int_{0}^{1} \frac{(1-t)^{p} (k+p+1)!}
  {\|x_{0}+t(x-x_{0})-y_{0}\|^{k+p+2}} \, dt,
\end{align*}
where the multinomial theorem and the rotational invariance were used
twice. The sum can be evaluated when the upper limit tends to
$\infty$; hence by uniform convergence,
\begin{align*}
  |\Efl| &\le C (p+1) \|x-x_{0}\|^{p+1} \int_{0}^{1} \frac{(1-t)^{p}}{\left(
      \|x_{0}+t(x-x_{0})-y_{0}\|-\|y-y_{0}\|\right)^{p+2}} \, dt \\
  &\le C (p+1) \frac{\|x-x_{0}\|^{p+1}}{\left(
      \|x_{0}-y_{0}\|-\|y-y_{0}\| \right)^{p+2}}
  \int_{0}^{1} \frac{(1-t)^{p}}{\left(
      1-t \frac{\|x-x_{0}\|}{\|x_{0}-y_{0}\|-\|y-y_{0}\|} \right)^{p+2}} \, dt.
\end{align*}
Using \eqref{eq:Gbound}, the triangle inequality, and \eqref{eq:ineq1}
we get (compare with \eqref{eq:relG1})
\begin{align}
  \label{eq:relG2}
  \frac{1}{\|x_{0}-y_{0}\|-\|y-y_{0}\|} 
  \le cG(x,y) \frac{\|x-y\|}{\|x_{0}-y_{0}\|-\|y-y_{0}\|} 
  \le cG(x,y) \frac{2}{1-\theta}.
\end{align}
For the integrand we use the same type of bound as before and finally get
\begin{align}
  \label{eq:Ef2l}
  |\Efl| &\le cC (p+1) \; \theta^{p+1} \;
  \frac{2}{(1-\theta)^{2}} G(x,y).
\end{align}

By summing the contributions \eqref{eq:Ef2f} and \eqref{eq:Ef2l} from
$N$ \emph{positive} potentials (cf.~\eqref{eq:Gbound}) we thus
conclude that the relative error for the $p$th order fast multipole
method under the \thetacriterion\ is bounded by a $\mbox{constant}
\times \theta^{p+1}/(1-\theta)^{2}$.

% *** (that \eqref{eq:Gbound} was used to get a relative error, can
% often get an absolute error estimate otherwise, but risking
% cancellation when the sign changes)

In the above derivation we are lead to the \thetacriterion\ by the
initial requirement that $\|y-y_{0}\|/\|x-y_{0}\| \le
\theta$. Clearly, the only natural way to obtain this is from the the
triangle inequality and the slightly stronger requirement that
$\|y-y_{0}\|/(\|x_{0}-y_{0}\|-\|x-x_{0}\|) \le \theta$
(cf.~\eqref{eq:ineq2}). Despite this clear reasoning we have not been
able to find our version of the \thetacriterion\ made explicit in the
literature. Presumably, this is due to the fact that quadratic meshes
are so popular.

We note also that the integral remainder term was consistently used
instead of the Lagrangian version (hinted at also in
\cite{fmmskeletons,ppfmm}). To see why, note that in bounding
e.g.~$\Eff$ above we would otherwise obtain terms of the form
$G_{0,\beta}(x,y_{0}+\xi(y-y_{0}))$ with $\xi \in [0,1]$ and $|\beta|
= p+1$. The only sensible bound now includes the factor
$(1-\theta)^{-(p+2)}$ suggesting that the convergence would
deteriorate as $\theta \to 1/2$. By contrast, the bounds
\eqref{eq:Ef2f} and \eqref{eq:Ef2l} are perfectly regular for
\emph{any} $0 < \theta < 1$.

It is to be stressed that we like to view the above analysis as a kind
of \emph{template}; the precise form of the assumptions and their
implications are all very clearly visible. The effect of changing them
can therefore readily be assessed.

% *** (more about the errors in GKZ?)

%**************************************************************************

\section{Implementation}
\label{sec:implementation}

\begin{figure}[htb]
  \centering
  \subfigure{\includegraphics{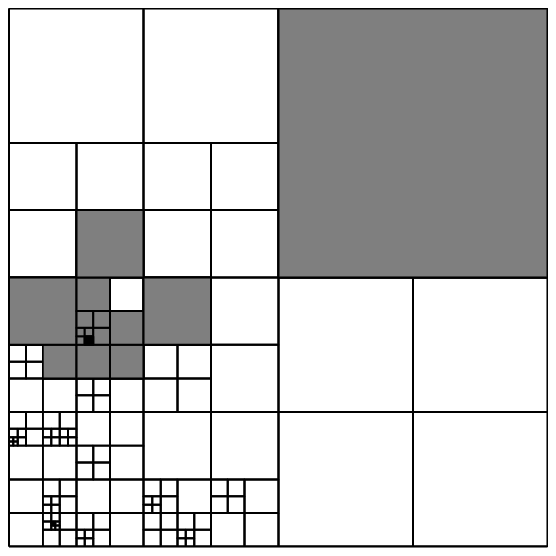}}
  \subfigure{\includegraphics{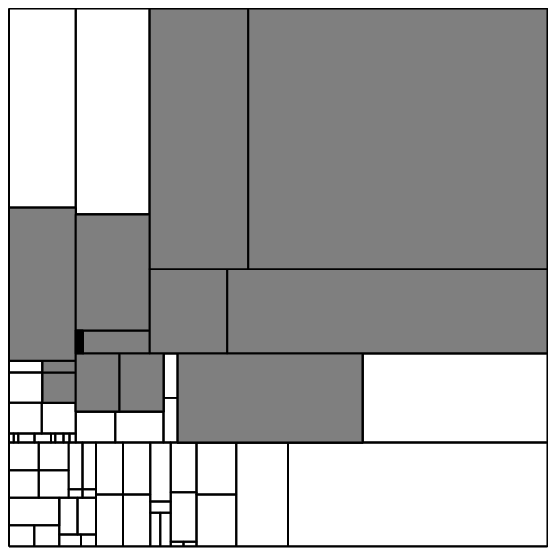}}
  \caption{Standard midpoint adaptivity \textit{(left)} vs.
    asymmetric adaptivity \textit{(right)} for the same set of point
    sources (not shown). In both cases the boxes that are strongly
    connected to the black box are shaded grey ($\theta = 1/2$). In
    the conventional approach, strongly connected boxes of different
    sizes imply cross-level communication in the multipole tree. By
    splitting the boxes at the median point source rather than at
    their geometrical midpoints, those boxes always belong to the same
    level.}
  \label{fig:boxes}
\end{figure}

The current paper on the \thetacriterion\ stems mainly from the fact
that allowing an asymmetric adaptive splitting when creating the FMM
mesh makes a convenient implementation possible. The ease of
implementation is mainly thanks to the fact that the associated
multipole tree is always balanced so that a static memory layout
becomes natural. As a result, neighboring boxes are always arranged
in the same level in the multipole tree, thereby facilitating the
communication between them (see Figure~\ref{fig:boxes} for an
illustration).

We now give a brief description of our current implementation. For
more detailed information the interested reader is kindly referred to
the freely available code itself (see Section~\ref{subsec:reproduce}).

The multipole mesh is constructed first and is obtained by recursively
splitting the source points at or near their median in each coordinate
direction. A very fast and well-known algorithm is available for doing
just this; this is the \emph{median-of-three} selection algorithm
which is most often used when implementing \textit{quicksort}
\cite[Chap.~9]{AlgorithmsC}. The source points themselves are
naturally stored in a \emph{pyramid} data-structure, a
\textit{$2^{d}$-tree} cut at a certain maximum level
\cite[Chap.~5.3]{helmholtzFMM}.

After the sources have been assigned a box at the finest level, the
\emph{connectivity} information is determined. At any level in the
tree the boxes are either decoupled or, respectively, strongly/weakly
coupled. For each box $b$, the strong connections $S(p)$ of its parent
box $p$ are examined; all children of a box in $S(p)$ that satisfy the
\thetacriterion\ with respect to $b$ become weakly coupled to $b$ ---
the rest remain strongly coupled. This simple rule together with the
fact that a box is always strongly connected to itself recursively
defines the connectivity for the whole tree. The resulting topology is
conveniently stored in sparse matrices with sparsity patterns known
\textit{before} each level is to be examined.

In the downward/upward phases the expansions are shifted as usual
according to parent/child relations. The critical far-to-local shift
follows the weak connectivity pattern and, at the finest level, the
remaining strong connections are evaluated directly. As suggested in
\cite{vertaxis2D,fmmBLAS} all shifts in the implementation tested
below relies on BLAS Level 3 routines with constant transition
matrices (using pre- and post-scaling according to the local
geometry).

A few structural optimizations are optionally possible. For instance,
a box weakly coupled to all children of a parent box can often
interact via the latter instead (whenever the \thetacriterion\ with
respect to that parent is true). This is visible in
Figure~\ref{fig:mesh} where some of the weakly coupled grey boxes are
noticeable large --- here the interactions are in fact managed via
parent boxes. A related idea at the \emph{finest} level is to
investigate strongly coupled boxes of very different radii, $r \ll R$,
say. If the \thetacriterion\ is true \emph{when the roles of $r$ and
  $R$ are exchanged}, then the sources in the larger box can be
directly converted into a near field expansion in the smaller box, and
the outgoing expansion from the smaller can simply be evaluated at
each point in the larger box. This optimization was suggested already
in \cite{AFMM}.

Since it is beneficial to admit as large set of boxes as possible into
the set of weak connections, it is natural to try to somehow locally
relax the \thetacriterion. Rather than uniformly enforcing a single
value of $\theta$ for example, \emph{different} $\theta_{i}$'s may
well be accepted provided that an estimate of the total error remains
bounded. For instance, if $n$ observed values of $\theta_{i}$ satisfy
$\sum_{i} \theta_{i}^{p+1}/n \le \theta^{p+1}$, then the relative
error estimate is still $\Ordo{\theta^{p+1}}$. Ideas along these lines
are discussed in \cite{fmmskeletons,rectangular_fmm}.

We conclude this section by very briefly discussing the algorithmic
complexity. It is known that there typically exists extremely
non-uniform distributions of $N$ source points that make most
tree-codes (even the adaptive ones) run in time proportional to
$\Ordo{N^{2}}$ \cite{FMMcomplexitynote}. Practical experience is
usually much better, although for certain applications, simpler
algorithms are occasionally preferred (see \cite{nbodycompare}, and
further the discussions in \cite[Chap.~8.7]{griebelMoldyn}, and
\cite[Chap.~6.6.3]{helmholtzFMM}).

In any case, the serial complexity for a 2D implementation using
asymmetric adaptivity can be estimated to be roughly proportional to
$\theta^{-2} p^{2} \cdot N$, since each of about $N$ boxes at the
finest level is to interact through cluster-to-cluster interactions
with on the order of $\theta^{-2}$ boxes at the same level, and since
each such shift requires on the order of $p^{2}$ operations. For a
given relative tolerance $\TOL$, the analysis in
Section~\ref{sec:analysis} implies $p \sim \log\TOL/\log\theta$, so
that the complexity is $\Ordo{\theta^{-2} \log^{-2}\theta \cdot
  N\log^{2} \TOL}$, thus indicating that the choice $\theta = \exp(-1)
\approx 0.368$ is nearly optimal. In practice, the best value of
$\theta$ as well as the optimal number of subdivisions is rather
dependent on the hardware and should be determined from experience. On
balance we have found that the convenient choice $\theta = 0.5$ and
subdividing the points until the number of sources per box is $\sim
20$ works very well in practice.

%**************************************************************************

\section{Experiments}
\label{sec:experiments}

As illustrations to the analysis in Section~\ref{sec:analysis} and in
order to highlight some of the benefits with the proposed adaptivity
we report some results from our two-dimensional implementation which
employs the classical complex polynomial/multipole representation
\cite{FMM,AFMM}.

\begin{SCfigure}[][htb]
  \includegraphics{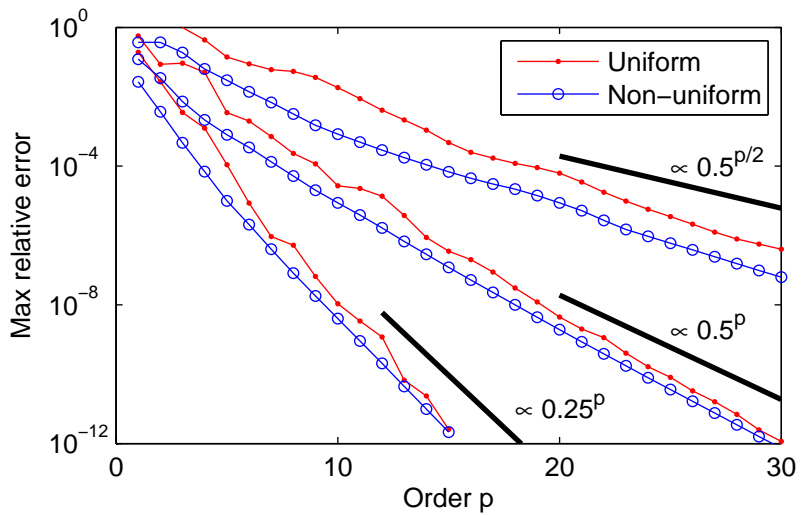}
  \caption{Relative error as a function of $p$ for $\theta \in
    \{0.25,0.5,\sqrt{0.5}\}$ and for two different distributions. A
    million point sources of total mass 1 are distributed inside the
    unit circle \textit{(i)} uniformly, and \textit{(ii)}
    non-uniformly with radial density $\propto 1/r$ (see
    Figure~\ref{fig:mesh} for the resulting mesh in the latter case).
    For convenience, the error is measured in a random sample ($M =
    1000$) of points.}
  \label{fig:errors}
\end{SCfigure}

Firstly, we investigated the sharpness of the \thetacriterion\ and the
accompanying error analysis. For this purpose, the complex-valued
field $G(z_{i},z_{j}) = -m_{j}/(z_{i}-z_{j})$ was used and we
calculated the total force in a system consisting of a million sources
(see Figure~\ref{fig:errors} for results and further details of this
simulation). Since $G$ is complex-valued the assumption
\eqref{eq:Gbound} is clearly violated. Cancellation effects therefore
implies that relative error estimates may well be impossible to
obtain. Nevertheless, since the points are distributed randomly and an
irregular mesh is used the effect is negligible in this case (although
somewhat more prominent for the uniform distribution, see
Figure~\ref{fig:errors}).

\begin{SCfigure}[][htb]
  \includegraphics{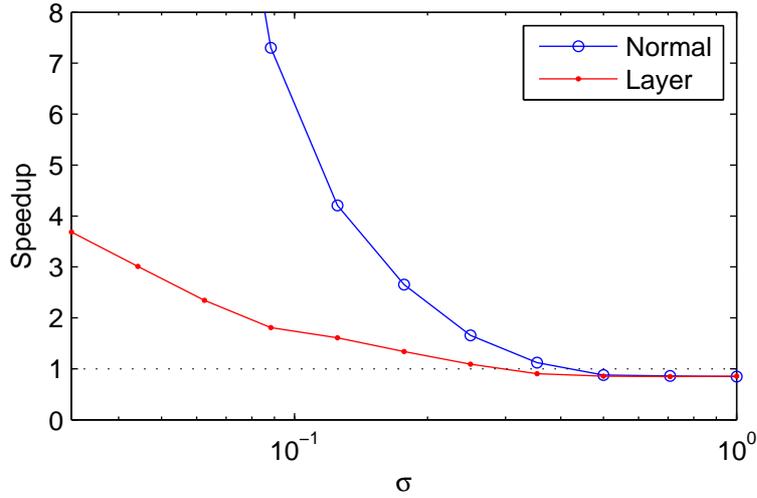}
  \caption{The benefits of adaptivity: speedup of the adaptive FMM
    versus an optimized uniform FMM. The harmonic potential is
    evaluated at 200,000 random points, chosen either from a normal
    distribution with variance $\sigma^{2}$ or from a combined `layer'
    distribution where the $x$-coordinate is uniformly distributed in
    $[0,1]$, and the $y$-coordinate again is
    $N(0,\sigma^{2})$-distributed. For ease of comparison, both
    distributions are forced by rejection to fit exactly within the
    unit square. As $\sigma \to 0$, stronger clustering around the
    origin and the $x$-axis occurs, respectively, and the benefits
    with adaptivity quickly show.}
  \label{fig:aavsu}
\end{SCfigure}

Secondly, we evaluated the efficiency of the proposed adaptivity at an
effective relative tolerance of about $10^{-6}$ (using $p = 20$
expansion coefficients). Since our code is developed from a highly
optimized uniform code outlined in \cite{vertaxis2D}, a reasonably
fair comparison is possible. Evidently, a square and uniformly
subdivided multipole mesh implies a completely regular access pattern
so that a highly efficient implementation using direct addressing
techniques is possible. It is therefore of some interest to estimate
under what circumstances adaptivity is actually beneficial. For this
purpose we measured the speedup achieved by the adaptive code when the
point sources were sampled from increasingly non-uniform data. The
results are displayed in Figure~\ref{fig:aavsu} and shows that only
for very uniform distributions of points is there a small gain
($\lesssim 15\%$) in using the uniform code.

\begin{SCfigure}[][htb]
  \includegraphics{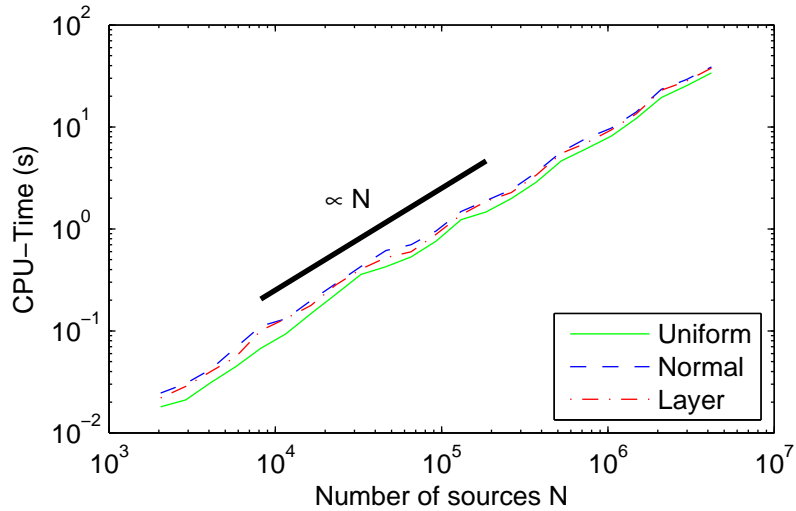}
  \caption{Robustness of the adaptive FMM for different distributions
    of sources. The CPU-time is measured for increasing number of
    sources and for three different distributions: a uniform
    distribution in the unit square, and the two distributions used
    also in Figure~\ref{fig:aavsu}.}
  \label{fig:robustness}
\end{SCfigure}

As a third and final experiment we assessed the performance of the
adaptive code for different point distributions. Using again a
relative tolerance of $10^{-6}$, we measured the CPU-time for
increasing number of points sampled from three very different
distributions. As shown in Figure~\ref{fig:robustness}, the code is
very robust indeed and scales well at least up to some 5 million point
sources on a single CPU.

\subsection{Reproducibility}
\label{subsec:reproduce}

Our two-dimensional implementation as described in
Section~\ref{sec:implementation} is available for download via the
author's
web-page\footnote{\url{http://user.it.uu.se/~stefane/freeware}}. The
code has been tested on several platforms, is fully documented, and
comes with a convenient Matlab mex-interface. Along with the code,
automatic Matlab-scripts that repeat the numerical experiments in
Section~\ref{sec:experiments} are also distributed. The results
presented here were all obtained with a 3.06 GHz Intel Core 2 Duo
processor with 4GB of internal memory running under Mac~OS~X~10.6.6.

%**************************************************************************

\section{Conclusions}

Asymmetric adaptive meshes offer convenient and efficient FMM
implementations. If cluster-to-cluster interactions are restricted to
boxes for which the \thetacriterion\ is true and if the kernel
satisfies the assumptions outlined in Section~\ref{sec:assumptions},
then the relative error can be bounded by a $\mbox{constant} \times
\theta^{p+1}/(1-\theta)^{2}$ for $\theta \in (0,1)$ and with $p$ the
number of expansion terms. The computational complexity of the
resulting 2D-algorithm can be expected to be about $\Ordo{\theta^{-2}
  \log^{-2}\theta \cdot N\log^{2} \TOL}$ for some target tolerance
$\TOL$. Not only is the actual performance of the algorithm
competitive with optimized uniform FMMs even for relatively uniform
data, but it is also robust for non-uniform distribution of points.
Ongoing work includes porting the code to manycore platforms.

%**************************************************************************

\section*{Acknowledgment}

This work was supported by the Swedish Research Council within the
FLOW and the UPMARC Linnaeus centers of Excellence.

%**************************************************************************

\bibliographystyle{model1b-num-names}
%\bibliography{../../../stefan}

\begin{thebibliography}{18}
\expandafter\ifx\csname natexlab\endcsname\relax\def\natexlab#1{#1}\fi
\providecommand{\bibinfo}[2]{#2}
\ifx\xfnm\relax \def\xfnm[#1]{\unskip,\space#1}\fi
%Type = Article
\bibitem[{Aluru(1996)}]{FMMcomplexitynote}
\bibinfo{author}{S.~Aluru}, \bibinfo{title}{Greengard's {$N$}-body algorithm is
  not order {$N$}}, \bibinfo{journal}{SIAM J.~Sci.~Comput.}
  \bibinfo{volume}{17} (\bibinfo{year}{1996}) \bibinfo{pages}{773--776}.
  \bibinfo{note}{\doi{10.1137/S1064827593272031}}.
%Type = Inproceedings
\bibitem[{Blelloch and Narlikar(1997)}]{nbodycompare}
\bibinfo{author}{G.~Blelloch}, \bibinfo{author}{G.~Narlikar}, \bibinfo{title}{A
  practical comparison of {$N$}-body algorithms}, in:
  \bibinfo{booktitle}{Parallel Algorithms}, volume~\bibinfo{volume}{30} of
  \textit{\bibinfo{series}{Series in Discrete Mathematics and Theoretical
  Computer Science}}.
%Type = Article
\bibitem[{Carrier et~al.(1988)Carrier, Greengard and Rokhlin}]{AFMM}
\bibinfo{author}{J.~Carrier}, \bibinfo{author}{L.~Greengard},
  \bibinfo{author}{V.~Rokhlin}, \bibinfo{title}{A fast adaptive multipole
  algorithm for particle simulations}, \bibinfo{journal}{SIAM
  J.~Sci.~Stat.~Comput.} \bibinfo{volume}{9} (\bibinfo{year}{1988})
  \bibinfo{pages}{669--686}. \bibinfo{note}{\doi{10.1137/0909044}}.
%Type = Article
\bibitem[{Cheng et~al.(1999)Cheng, Greengard and Rokhlin}]{new_AFMM}
\bibinfo{author}{H.~Cheng}, \bibinfo{author}{L.~Greengard},
  \bibinfo{author}{V.~Rokhlin}, \bibinfo{title}{A fast adaptive multipole
  algorithm in three dimensions}, \bibinfo{journal}{J.~Comput.~Phys.}
  \bibinfo{volume}{155} (\bibinfo{year}{1999}) \bibinfo{pages}{468--498}.
  \bibinfo{note}{\doi{10.1006/jcph.1999.6355}}.
%Type = Article
\bibitem[{Coulaud et~al.(2010)Coulaud, Fortin and Roman}]{fmmBLAS}
\bibinfo{author}{O.~Coulaud}, \bibinfo{author}{P.~Fortin},
  \bibinfo{author}{J.~Roman}, \bibinfo{title}{High performance {BLAS}
  formulation of the adaptive fast multipole method},
  \bibinfo{journal}{Math.~Comput.~Modelling} \bibinfo{volume}{51}
  (\bibinfo{year}{2010}) \bibinfo{pages}{177--188}.
  \bibinfo{note}{\doi{10.1016/j.mcm.2009.08.039}}.
%Type = Article
\bibitem[{Deglaire et~al.(2009)Deglaire, Engblom, {\AA}gren and
  Bernhoff}]{vertaxis2D}
\bibinfo{author}{P.~Deglaire}, \bibinfo{author}{S.~Engblom},
  \bibinfo{author}{O.~{\AA}gren}, \bibinfo{author}{H.~Bernhoff},
  \bibinfo{title}{Analytical solutions for a single blade in vertical axis
  turbine motion in two-dimensions}, \bibinfo{journal}{Eur.~J.~Mech.~B Fluids}
  \bibinfo{volume}{28} (\bibinfo{year}{2009}) \bibinfo{pages}{506--520}.
  \bibinfo{note}{\doi{10.1016/j.euromechflu.2008.11.004}}.
%Type = Article
\bibitem[{Greengard and Rokhlin(1987)}]{FMM}
\bibinfo{author}{L.~Greengard}, \bibinfo{author}{V.~Rokhlin}, \bibinfo{title}{A
  fast algorithm for particle simulations}, \bibinfo{journal}{J.~Comput.~Phys.}
  \bibinfo{volume}{73} (\bibinfo{year}{1987}) \bibinfo{pages}{325--348}.
  \bibinfo{note}{\doi{10.1016/0021-9991(87)90140-9}}.
%Type = Book
\bibitem[{Griebel et~al.(2007)Griebel, Knapek and Zumbusch}]{griebelMoldyn}
\bibinfo{author}{M.~Griebel}, \bibinfo{author}{S.~Knapek},
  \bibinfo{author}{G.~Zumbusch}, \bibinfo{title}{Numerical Simulation in
  Molecular Dynamics}, volume~\bibinfo{volume}{5} of
  \textit{\bibinfo{series}{Texts in Computational Science and Engineering}},
  \bibinfo{publisher}{Springer Verlag}, \bibinfo{address}{Berlin},
  \bibinfo{year}{2007}.
%Type = Book
\bibitem[{Gumerov and Duraiswami(2004)}]{helmholtzFMM}
\bibinfo{author}{N.A. Gumerov}, \bibinfo{author}{R.~Duraiswami},
  \bibinfo{title}{Fast multipole methods for the {H}elmholtz equation in three
  dimensions}, Elsevier Series in Electromagnetism,
  \bibinfo{publisher}{Elsevier}, \bibinfo{address}{Oxford},
  \bibinfo{year}{2004}.
%Type = Article
\bibitem[{Gumerov and Duraiswami(2008)}]{fmmgpu}
\bibinfo{author}{N.A. Gumerov}, \bibinfo{author}{R.~Duraiswami},
  \bibinfo{title}{Fast multipole methods on graphics processors},
  \bibinfo{journal}{J.~Comput.~Phys.} \bibinfo{volume}{227}
  (\bibinfo{year}{2008}) \bibinfo{pages}{8290--8313}.
  \bibinfo{note}{\doi{10.1016/j.jcp.2008.05.023}}.
%Type = Article
\bibitem[{Petersen et~al.(1995)Petersen, Smith and Soelvason}]{fmmerr}
\bibinfo{author}{H.G. Petersen}, \bibinfo{author}{E.R. Smith},
  \bibinfo{author}{D.~Soelvason}, \bibinfo{title}{Error estimates for the fast
  multipole method. {II}. {T}he three-dimensional case},
  \bibinfo{journal}{Proc.~Math.~Phys.~Sci.} \bibinfo{volume}{448}
  (\bibinfo{year}{1995}) \bibinfo{pages}{401--418}.
%Type = Article
\bibitem[{Salmon and Warren(1994)}]{fmmskeletons}
\bibinfo{author}{J.K. Salmon}, \bibinfo{author}{M.S. Warren},
  \bibinfo{title}{Skeletons from the treecode closet},
  \bibinfo{journal}{J.~Comput.~Phys.} \bibinfo{volume}{111}
  (\bibinfo{year}{1994}) \bibinfo{pages}{136--155}.
  \bibinfo{note}{\doi{10.1006/jcph.1994.1050}}.
%Type = Book
\bibitem[{Sedgewick(1990)}]{AlgorithmsC}
\bibinfo{author}{R.~Sedgewick}, \bibinfo{title}{Algorithms in {C}},
  Addison-Wesley Series in Computer Science,
  \bibinfo{publisher}{Addison-Wesley}, \bibinfo{address}{Reading, MA},
  \bibinfo{year}{1990}.
%Type = Article
\bibitem[{Shanker and Huang(2007)}]{ace}
\bibinfo{author}{B.~Shanker}, \bibinfo{author}{H.~Huang},
  \bibinfo{title}{Accelerated {C}artesian expansions - a fast method for
  computing of potentials of the form ${R}^{-\nu}$ for all real $\nu$},
  \bibinfo{journal}{J.~Comput.~Phys.} \bibinfo{volume}{226}
  (\bibinfo{year}{2007}) \bibinfo{pages}{732--753}.
  \bibinfo{note}{\doi{10.1016/j.jcp.2007.04.033}}.
%Type = Article
\bibitem[{Sundar et~al.(2008)Sundar, Sampath and Biros}]{balancedtrees}
\bibinfo{author}{H.~Sundar}, \bibinfo{author}{R.S. Sampath},
  \bibinfo{author}{G.~Biros}, \bibinfo{title}{Bottom-up construction and 2:1
  balance refinement of linear octrees in parallel}, \bibinfo{journal}{SIAM
  J.~Sci.~Comput.} \bibinfo{volume}{30} (\bibinfo{year}{2008})
  \bibinfo{pages}{2675--2708}. \bibinfo{note}{\doi{10.1137/070681727}}.
%Type = Inproceedings
\bibitem[{Vikram et~al.(2009)Vikram, Baczewzki, Shanker and Aluru}]{pace}
\bibinfo{author}{M.~Vikram}, \bibinfo{author}{A.~Baczewzki},
  \bibinfo{author}{B.~Shanker}, \bibinfo{author}{S.~Aluru},
  \bibinfo{title}{Parallel accelerated {C}artesian expansions for particle
  dynamics simulations}, in: \bibinfo{booktitle}{Proceedings of the 2009 IEEE
  International Parallel and Distributed Processing Symposium}, pp.
  \bibinfo{pages}{1--11}. \bibinfo{note}{\doi{10.1109/IPDPS.2009.5161038}}.
%Type = Article
\bibitem[{Warren and Salmon(1995)}]{ppfmm}
\bibinfo{author}{M.S. Warren}, \bibinfo{author}{J.K. Salmon}, \bibinfo{title}{A
  portable parallel particle program}, \bibinfo{journal}{Comput.~Phys.~Commun.}
  \bibinfo{volume}{87} (\bibinfo{year}{1995}) \bibinfo{pages}{266--290}.
  \bibinfo{note}{\doi{10.1016/0010-4655(94)00177-4}}.
%Type = Article
\bibitem[{Zhang and Tanaka(2007)}]{rectangular_fmm}
\bibinfo{author}{J.M. Zhang}, \bibinfo{author}{M.~Tanaka},
  \bibinfo{title}{Adaptive spatial decomposition in fast multipole method},
  \bibinfo{journal}{J.~Comput.~Phys.} \bibinfo{volume}{226}
  (\bibinfo{year}{2007}) \bibinfo{pages}{17--28}.
  \bibinfo{note}{\doi{10.1016/j.jcp.2007.03.032}}.

\end{thebibliography}

\providecommand{\noopsort}[1]{} \providecommand{\doi}[1]{\texttt{doi:#1}}
  \providecommand{\available}[1]{Available at \texttt{#1}}
  \providecommand{\availablet}[2]{Available at \texttt{#2}}

\end{document}